\documentclass[11pt]{article}
\usepackage{amsfonts}
\usepackage{latexsym}
\usepackage{amsmath}
\usepackage{amssymb}
\usepackage{color}
\usepackage{amsthm}
\usepackage{fullpage}
\usepackage{enumerate}
\newtheorem{theorem}{Theorem}[section]
\newtheorem{lemma}[theorem]{Lemma}

\newtheorem{remark}[theorem]{Remark}
\newtheorem{conjecture}[theorem]{Conjecture}
\theoremstyle{definition}

\newtheorem{thmy}{Theorem}
\newenvironment{oldtheorem}{\stepcounter{thm}\begin{thmy}}{\end{thmy}}
\newtheorem*{note*}{Note}

\makeatletter

\newcommand{\sign}{\mbox{sign}}
\newcommand{\R}{\mathbb{R}}
\newcommand{\Sp}{\mathbb{S}}
\newcommand{\N}{\mathbb{N}}
\newcommand{\RR}{\mathcal R}
\newcommand{\Q}{\mathcal Q}

\def\endproof{\begin{flushright}
$ \Box $ \\
\end{flushright}}
\newcommand{\beginproof}{\noindent{\bf Proof: }}

\makeatother

\makeatletter
\newcommand{\subjclass}[2][1991]{%
  \let\@oldtitle\@title%
  \gdef\@title{\@oldtitle\footnotetext{#1 \emph{Mathematics subject classification.} #2}}%
}
\newcommand{\keywords}[1]{%
  \let\@@oldtitle\@title%
  \gdef\@title{\@@oldtitle\footnotetext{\emph{Key words and phrases.} #1.}}%
}
\makeatother

\begin{document}

\title{Iterations of the projection body operator and a remark on Petty's conjectured projection inequality}

\author{C. Saroglou and A. Zvavitch\thanks{Supported in part by the U.S. National Science Foundation Grant DMS-1101636 and by the  Simons Foundation.}}

\subjclass[2010]{52A20, 53A15, 52A39. }
 \keywords{Convex bodies;   Projection Bodies; Petty's projection inequality.}

\date{\today}
\maketitle
\abstract{We prove that if a convex body has absolutely continuous surface area measure, whose density is sufficiently close to the constant, then the sequence $\{\Pi^mK\}$ of convex bodies converges to the ball with respect to the Banach-Mazur distance, as $m\rightarrow\infty$. Here, $\Pi$ denotes the projection body operator.  Our result allows us to show that the ellipsoid is a local  
solution to the conjectured inequality of Petty and to improve a related inequality of Lutwak.}
\section{Introduction}

As usual,   we denote by $\langle x,y\rangle$  the inner product of two vectors $x, y \in \R^d$ and by $|x|$ the length of a vector $x \in \R^d$. A {\em body} is   a compact  set with nonempty interior. For a body $K$ which is star-shaped with respect to the origin its {\em radial function}  is defined by
$$
\rho_K(u) = \max \{t \ge 0 : t u \in K\} \quad\mbox{for every } u\in \Sp^{d-1}.
$$
A body $K$ is called a {\em star body} if it is star-shaped at the origin and its radial function $\rho_K$ is positive and continuous.
 We say that a set $K$ is {\it symmetric} if it is centrally symmetric with center at the origin, i.e. for every $x \in K$ we get $-x \in K$.
We write $|A|$ for the $k$-dimensional Lebesgue measure (volume)  of a measurable set $A \subset \R^d$, where $k=1,\dots, d$  is the dimension of the minimal flat containing $A$. We recall also the definition of the {\it Banach-Mazur distance} between  symmetric  bodies $K$ and $L$
  $$d_{BM}(K,L) = \inf \{t\ge1 : L\subset TK\subset tL, \mbox{ for some } T\in{\rm GL}(d)\}.$$
The above distance turns out to be extremely useful when  one would like to study questions invariant under the linear transformation.  For a convex body $K$, the  support function $h_K(\xi):\Sp^{d-1} \to \R_+$ is defined as
$$ 
h_K(\xi)=\max\limits_{x\in K}  \langle x,\xi \rangle.
$$
We refer to \cite{G, Gr, Ko, KoY, Sch1} for definitions and properties of star-shaped and convex bodies and corresponding functionals.

The projection body of a convex body $K$ in $\mathbb{R}^d$ is defined as the body with support function
$$h_{\Pi K}(x)=|K|x^{\bot}|,  \mbox{ for all   }  x\in \Sp^{d-1},$$
where $K|x^{\bot}$ denotes the orthogonal projection of $K$ onto the subspace $x^{\bot}=\{y\in \R^d:  \langle x,y\rangle =0 \}$.  The direct application of Cauchy projection formula gives us
$$h_{\Pi K}(x)=\frac{1}{2}\int_{\Sp^{d-1}}|\langle x,y\rangle |dS_K(y),\ x\in \Sp^{d-1},$$
where $S_K$ is the surface area measure of $K$, viewed as a measure on $\Sp^{d-1}$. When $S_K$ is absolutely continuous (with respect to the Lebesgue measure on the sphere), its density $f_K$ is called curvature function of $K$.

It is very interesting to study the iterations of projection body operator. It is trivial to see that the projection body of Euclidean ball $B_2^d$ is again, up to a dilation,   $B_2^d$, moreover the same property is true for a 
unit cube $B_\infty^d$.  Weil \cite{W1} proved that if $K$ is a polytope in $\mathbb{R}^d$, then $\Pi^2K$ is homothetic to $K$ if and only if $K$ is a linear image of cartesian products of planar symmetric polygons or line segments.  But no other description of  fixed points of projection body operator is known as well as no much known about possible convergence of the sequence $\Pi^mK$. Clearly, Weil's result tell us that   one cannot expect in general that $\Pi^mK\rightarrow B^d_2$, with respect to the Banach-Mazur distance. It seems more plausible, however, that $\Pi^mK\rightarrow B^d_2$, if $K$ has absolutely continuous surface area measure and $d\geq 3$ (for $d=2$, if $K$ is symmetric, then $\Pi^2 K=4K$). We refer to \cite{G, Sch1} for more information about this problem. We are going to show the following:
\begin{theorem}\label{main-theorem}
Let $d\geq 3$. There exists an $\varepsilon_d>0$ with the following property: For any convex body $K$ in $\mathbb{R}^d$, with absolutely continuous surface area measure and   the curvature function  $f_K$ satisfying  $\|f_{TK}-1\|_{\infty}<\varepsilon_d$, for some $T\in GL(d)$, we have
$\Pi^mK\rightarrow B_2^d$, in the sense of the Banach-Mazur distance.
\end{theorem}

The idea of the above theorem follows from the study of the properties of intersection body operator done in \cite{FNRZ}. The  {\em intersection body} $IK$ of a star body $K$ was  defined by  Lutwak in \cite{Lu1}  using the radial function of the  body $IK$:
 $$\rho_{IK}(u)= |K\cap u^\bot| , \quad \mbox{ for   } u \in \Sp^{d-1}.$$
Again it is trivial to see that $IB_2^d$ is a dilate of $B_2^d$ and $I^2K= 4K$ for symmetric $K\subset \R^2$,  but no much information is known about other fixed points of $I$  (see \cite{Lu4}). In  \cite{FNRZ} authors shown that $B_2^d$ is a local attractor:
\begin{oldtheorem}\label{oldthm1}Let $d\geq 3$. There exists an $\varepsilon_d>0$ with the following property: For any star body $K$ in $\mathbb{R}^d$, which satisfies $\|\rho_{TK}-1\|_{\infty}<\varepsilon_d$, for some $T\in GL(d)$ (in other words, $K$ is close, in Banach-Mazur distance,  to  $B_2^d$ ), we have
$I^mK\rightarrow B_2^d$, in the sense of the Banach-Mazur distance.
\end{oldtheorem}

Another reason to consider   Theorem  \ref{main-theorem} is that it can be applied to study of Petty's conjectured inequality. Indeed, it was shown by Petty \cite{PE2} that the quantity $P(K):=|\Pi K|/|K|^{d-1}$ is affine invariant. Petty \cite{PE} also conjectured the following:
\begin{conjecture}
Let $d\geq 3$. The affine invariant $P(K)$ is minimal if and only if $K$ is an ellipsoid.
\end{conjecture}
The restriction $d\geq 3$ is because in the plane it is well known that $|\Pi K|\geq 4|K|$, with equality if and only if $K$ is symmetric (see, for example, 
\cite{Sa1}). Petty's conjecture, if true, would be a very strong inequality, as it would imply a number of important isoperimetric inequalities, such as the classical isoperimetric inequality, the Petty projection inequality (see \cite{Zh} for a functional form) and the affine isoperimetric inequality. 


 Very little seem to be known about the conjecture of Petty. For instance, as shown in \cite{Sa1}, Steiner symmetrization fails for this problem. A useful fact, that will be used subsequently-established by Schneider, is that 
\begin{equation}\label{eq:schn}
P(K)\geq P(\Pi K),
\end{equation}
 with equality if and only if $K$ is homothetic to $\Pi^2 K$ (see \cite{Sch3} or \cite[pp 570]{Sch1}). In particular, it follows that every solution to the Petty problem must be a zonoid (a body which is a limit of Minkowski sum of segments).

Although bodies with minimal surface area significantly larger than the surface area of the ball (of the same volume) satisfy the Petty conjecture (see e.g. \cite{GP}), no natural class of convex bodies was known to satisfy the Petty conjecture (natural class means that is connected with respect to the Banach-Mazur distance and contains the ball). Below, we have a result towards this direction.

\begin{theorem}\label{corollary 1}
Let $d\geq 3$. There exists an $\varepsilon_d>0$ with the following property: For any non-ellipsoidal convex body $K$ in $\mathbb{R}^d$, which has absolutely continuous surface area measure and satisfies $\|f_{TK}-1\|_{\infty}<\varepsilon_d$, for some $T\in GL(d)$, we have
$P(K)> P(B_2^d)$.
\end{theorem}
\beginproof Let $K$ be a convex body that satisfies the assumptions of Theorem \ref{corollary 1}. The functional $P$ is continuous with respect to the Banach-Mazur distance (this follows immediately from continuity of volume measure and the continuity of the projection body operation, see \cite{Sch1} for more details),  so by Theorem \ref{main-theorem},
$P(\Pi^mK)\rightarrow P(B_2^d)$. Also, the sequence $P(\Pi^mK)$ is non-increasing, thus $$P(K)\geq P(\Pi K)\geq  \lim_{m\rightarrow\infty}P(\Pi^mK)=P(B_2^d),$$ with equality in the first inequality if and only if $K$ is homothetic to $\Pi^2K$ (see inequality (\ref{eq:schn}) from above). But if $K$ is homothetic to $\Pi^2K$, then trivially $\Pi^{2m}K\rightarrow K$, with respect to the Banach-Mazur distance, which shows that $K$ is an ellipsoid. This is a contradiction and our assertion follows. 
\endproof
One might hope that Theorem \ref{main-theorem} can more generally show that for every convex body $K$, which is sufficiently close to the ellipsoid with respect to the Banach-Mazur distance, we have $P(K)\geq P(B_2^d)$. The difficulty here is that even if we assume that $K$ is smooth enough, the fact that $K$ is close to $B_2^d$ does not guarantee that $f_K$ is close to $f_{B_2^d}$. Indeed, the curvature function $f_K$ can be written as a determinant involving second derivatives of $h_K$. Thus to guarantee that $f_K$ is close to $f_{B_2^n}$ we need to have a restriction on second derivatives  of $h_K$, which is an almost equivalent statement to the fact that $f_K$ is close to $f_{B_2^d}$. Unfortunately,  we also do not see how direct approximation methods may help. However, the next remark explains how one can restate Theorem \ref{corollary 1} involving the Banach-Mazur distance. 
\begin{remark}
As we discussed above, it is enough to solve  the conjecture of Petty  in the class of zonoids. Moreover,  to prove the inequality, it is enough to consider  only zonoids with absolutely continuous surface area measure. This class corresponds to the class of centroid bodies of star bodies.  Therefore, Petty's conjectured inequality (without equality cases) is equivalent to the following statement: the affine invariant $S(K):=|\Pi \Gamma K|/|\Gamma K|^{d-1}$ is minimal among all star bodies $K$ if $K$ is an ellipsoid.
Here, $\Gamma K $ denotes the centroid body of $K$ i.e. the convex body whose support function is given by
$h_{\Gamma K}(x)=\int_K|\langle x,y\rangle|dy$, $x\in\mathbb{R}^d$.
Theorem \ref{corollary 1} says that if the star body $K$ is close enough (in the sense of the Banach-Mazur distance) to the ball, then $S(K)\geq S(B_2^d)$. 
\end{remark}

Denote by $W_i$ the $i$-th quermassintegral functional in $\mathbb{R}^d$, $i=0,1,\dots,d-1$, which is the mixed volume of $d-i$ copies of a convex body $K$ with $i$ copies of $B_2^d$ (see \cite{Sch1} for exact definitions and properties). Recall the Aleksandrov-Fenchel inequalities for quermassintegrals:
\begin{equation}\label{A-F-inequalities}
W_{i+1}^{d-i}(K)\geq\omega_dW_i^{d-i-1}(K),\qquad i=0,\dots,d-2,
\end{equation}
where $K$ is any convex body and $\omega_d=|B_2^d|$ (see \cite{Sch1}). Lutwak \cite{Lut2} proved that if Petty's conjecture was proven to be true, then a family of inequalities that are stronger than (\ref{A-F-inequalities}) would have been established. These conjectured inequalities involve the notion of the $i$-th projection body $\Pi_iK$ of $K$, whose support function is given by:
$$h_{\Pi_iK}(u)=W_{i|u^{\bot}}(K|u^{\bot}), \qquad i=0,\dots, d-2,$$
where $W_{i |u^{\bot}}$ stands for the $i$-th quermassintegral in $u^{\bot}$. Note that
$\Pi K=\Pi_0K$.
Actually, Lutwak  in \cite{Lut2}  established a certain member of this family of inequalities:
\begin{equation}\label{lut-ineq}
W_{d-2}(\Pi_{d-2}K)\geq \omega_{d-1}^2W_{d-2}(K), 
\end{equation}
where $d\geq 3$, with equality if and only if $K$ is a ball.
To see that (\ref{lut-ineq}) is stronger than (\ref{A-F-inequalities}) (in the sense that it interpolates  (\ref{A-F-inequalities})), for $i=d-1$, note that since $W_{d-1}(K)$ is proportional to the mean width of $K$ (see
also \cite{Lut2}), we  get:
\begin{equation}
 W_{d-1}(\Pi_{d-2}K)=\omega_{d-1}W_{d-2}(K),\label{identity}
\end{equation}
Thus by (\ref{A-F-inequalities}) we obtain:
$$\frac{\omega_d}{\omega_{d-1}^2}W_{d-2}(\Pi_{d-2}K)\leq W_{d-1}^2(K),$$
with equality if and only if $\Pi_{d-2}K$ is a ball.

Theorem \ref{corollary 1}  allows us to prove a stronger version of Lutwak's inequality:
\begin{theorem}\label{corollary 2}
Let $K$ be a convex body in $\mathbb{R}^d$, $d\geq 3$. Then,
$$W_{d-2}(\Pi_{d-2}K)\geq \frac{d(d-2)\omega_{d-1}^2}{(d-1)^2\omega_d}W^2_{d-1}(K)+\frac{\omega_{d-1}^2}{(d-1)^2}W_{d-2}(K).$$
This inequality is sharp for the ball. Moreover, if $K$ is not centrally symmetric, then the inequality is strict.
\end{theorem}
As noted in \cite{Lut2}, since $W_0(K)=|K|$ and $\Pi_0K=\Pi K$, in the plane (\ref{lut-ineq}) becomes $|\Pi K|\geq 4|K|$, which as mentioned earlier is known to hold and, also, the inequality is strict if $K$ is not centrally symmetric. Therefore, Theorem \ref{corollary 2} is trivially true in dimension 2 as well. 

The paper is organized as follows: in Section 2 we use the method of Spherical Harmonics together with basic lemmata from \cite{FNRZ} to present a proof of Theorem 1.1 and  Section 3 is dedicated to  the proof of Theorem \ref{corollary 2}.
\\

\noindent {\bf Acknowledgment}: We are indebted to Fedor Nazarov  who suggested to us the direct link between $\Q$ transform and the Radon transform (equation (\ref{observation}) below) and Dmitry Ryabogin for many
valuable discussions.

\section{Iterations of the projection body operator}
The goal of this section is to establish Theorem \ref{main-theorem}.
Let us first fix some notation. If $f:\Sp^{d-1}\rightarrow \mathbb{R}$, denote by $\RR(f)$ the normalized Radon transform of $f$ (i.e. $\RR(1)=1$). In other words, 
$$
\RR(f)(x)=\frac{1}{|\Sp^{d-2}|}\int_{\Sp^{d-1}\cap x^{\bot}}fdy,\qquad x\in \Sp^{d-1}.
$$
Consider functions  $g_1,\dots,g_{d-1}:\Sp^{d-1}\rightarrow\mathbb{R}$ and define the function $\Q(g_1,\dots,g_{d-1}): \Sp^{d-1}\mapsto \R$ as :
$$\Q(g_1,\dots,g_{d-1})(x)=c_d\int\limits_{\Sp^{d-1}\cap x^{\bot}}\dots\int\limits_{\Sp^{d-1}\cap x^{\bot}}\det(x_1,\dots,x_{d-1})^2g_1(x_1)\dots g_{d-1}(x_{d-1})dx_1\dots dx_{d-1}.$$
We abbreviate $\Q(f):=\Q(f,\dots,f)$ and 
$$
\Q(g_1[k],g_2[d-1-k])=\Q(\underset{k\ {\rm times}}{\underbrace{g_1,\dots,g_1}},\underset{d-1-k\ {\rm times}}{\underbrace{g_2,\dots,g_2}}).$$
 The constant $c_d$ is chosen so that $\Q(1)=1$.

Denote, also, by $H_k^f$ the $k$-th spherical harmonic of $f$, $k\in\mathbb{N}$. It is well known (see e.g. \cite{Gr},  Lemma 3.4.7, page 103)
that $\RR(H_k^f)=(-1)^{k/2}v_{k,d}H_k^f$, if $k$ is even and $\RR(H_k^f)=0$ if $k$ is odd, where
$$v_{k,d}=\frac{1\cdot 3\dots (k-1)}{(d-1)(d+1)\dots (d+k-3)}.$$
Note that $H_0^f=|\Sp^{d-1}|^{-1}\int_{\Sp^{d-1}}fdx$ and that $\|f-H_0^f\|_{\infty}\leq 2\|f\|_{\infty}$.

It will be essential for us to consider  different actions of a linear operator  $T\in GL(d)$ on a function $f:\Sp^{d-1} \to \R$. We  define  $\widehat{Tf}=f(Tx/|Tx|)|Tx|^{-1}$ and $Tf=f(Tx/|Tx|)|Tx|^{-(d+1)}$.

Let $K$ be a convex body with  a curvature function $f_K$ it was proved  by Lutwak \cite{Lut3}, that if $T\in SL(d)$, then $f_{TK}=T^{\ast}f_K$. On the other hand, for $t>0$, we have $f_{tK}=t^{d-1}f_K$.
Therefore, for $T\in GL(d)$, since $T/|\textnormal{det}T|^{1/d}\in SL(d)$, it follows easily that $$f_{T^{\ast}K}=|\textnormal{det}T|^{(d-1)/d}Tf_K.$$ Moreover, it was shown by Weil \cite{W2} that projection body of $K$ again has absolutely continuous surface area measure and its density is proportional to
$\Q(f_K)$. For this section, it will be convenient to normalize $\Pi K$, so that $\Q(f_K)=f_{\Pi K}$.

The proof of Theorem \ref{main-theorem} will be done through a modification of the proof of Theorem A. All constants that appear in this section will depend on the dimension $d$ only. For $a\geq 0$, the following norm was introduced in \cite{FNRZ}: 
\begin{align*}
\|f\|_{{\cal U}_a}&=\\
\inf\{&M>0:\|f\|_{\infty}<M\textnormal{ and  } \forall n\in \N,\ \exists\ p_n \textnormal{ polynomial of degree } n, \textnormal{ s.t. } \|f-p_n\|_2<Mn^{-a}\}.
\end{align*}
If $\|f\|_{{\cal U}_a}<\infty$, we say that $f$ belongs to the class ${\cal U}_a$.
It was shown in \cite[proof of Lemma 6, Step 1]{FNRZ} that there exist constants $a, L, C>0$, with the following properties: 
\begin{itemize}
\item[(i)] If $\varphi:\Sp^{d-1}\rightarrow \mathbb{R}$ is a function with $\|\varphi\|_2<\varepsilon$, for some $0<\varepsilon<1$ and $\|\varphi\|_{{\cal U}_a}<1/L$, then $\|\varphi\|_{\infty}<C\varepsilon^{\frac{4}{d+3}}$. \item[(ii)] If $f=1+\varphi$ and $T\in GL(d)$ with $T=I+Q$,  for some $Q\in GL(d)$,  then
$$\widehat{Tf}=1-\langle Qx,x\rangle+\varphi(x)+z(x),$$
where $\|z\|_{\infty}=O(|Q|\varepsilon^{\frac{2}{d+3}}+|Q|^2)$,
\end{itemize}
where by $|Q|$ we denote the operator norm of $Q$. Thus, since
$$|Tx|^{-d}=1-d\langle Qx,x\rangle+O(|Q|^2),$$
we get
\begin{eqnarray*}
Tf(x)&=&1-(d+1)\langle Qx,x\rangle+\varphi(x)+d\langle Qx,x\rangle\varphi(x)+O(|Q|\varepsilon^{\frac{2}{d+3}}+|Q|^2)\\
&=&1-(d+1)\langle Qx,x\rangle+\varphi(x)+O(|Q|\varepsilon^{\frac{2}{d+3}}+|Q|^2).
\end{eqnarray*}
Choose $Q$ such that $-(d+1)\langle Qx,x\rangle=H_2^{\varphi}(x)$. Then, $\|H_2^{\varphi}(x)\|_2\leq \|\varphi\|_2<\varepsilon$. 
Write $\langle Qx,x\rangle=\sum_{i,j=1}^dq_{ij}x_ix_j$. Then, since $\int_{\Sp^{d-1}}x_1x_2x_3x_4dx=\int_{\Sp^{d-1}}x_1^2x_2x_3dx=0$, we have:
$$\|\langle Qx,x\rangle\|_2^2=\sum_{i\neq j}q_{ij}^2\|x_ix_j\|_2^2+\sum_{i=1}^dq_{ii}^2\|x_i^2\|_2^2=l_1\sum_{i\neq j}q_{ij}^2+l_2\sum_{i=1}^dq_{ii}^2\geq \min\{l_1,l_2\}q_{ij}^2,$$
for all $i,j$, where $l_1,l_2$ are positive constants that depend only on the dimension. This easily implies that $|Q|=O(\varepsilon)$. Therefore,
$Tf=1+\psi$, with
$\psi=\varphi-H_2^{\varphi}+O(\varepsilon^{\frac{d+5}{d+3}})$. Thus, $\|\psi\|_2\leq \|\varphi\|_2+O(\varepsilon^{\frac{d+5}{d+3}})\leq \varepsilon+C'\varepsilon^{\frac{d+5}{d+3}}$ and $\|\psi\|_{\infty}\leq \|\varphi\|_{\infty}+\|H_2^{\varphi}\|_{\infty}+O(\varepsilon^{\frac{d+5}{d+3}})\leq C''\varepsilon^{\frac{4}{d+3}}$. Finally, since $\varphi-H_2^{\varphi}$ has no second degree spherical harmonic in its expansion, we have $\|H_2^{\psi}\|_2\leq C'''\varepsilon^{\frac{d+5}{d+3}}$. Therefore, we have the following:
\begin{lemma}\label{lemma 1}
 There exist constants $a, L, C>0$, with the following property: Whenever a function $\varphi:\Sp^{d-1}\rightarrow \mathbb{R}$ satisfies $\|\varphi\|_2<\varepsilon$, for some $0<\varepsilon<1$ and $\|\varphi\|_{{\cal U}_a}<1/L$, then (i) $\|\varphi\|_{\infty}<C\varepsilon^{\frac{4}{d+3}}$ and (ii) there exists $T\in GL(d)$, such that if we set $f=1+\varphi$, then $Tf=1+\psi$, with  $\|\psi\|_2\leq \varepsilon+C\varepsilon^{\frac{d+5}{d+3}}$, $\|\psi\|_{\infty} \le C\varepsilon^{\frac{4}{d+3}}$ and
$\|H_2^{\psi}\|_2\leq C\varepsilon^{\frac{d+5}{d+3}}$.
\end{lemma}
We will also need two more lemmas from \cite{FNRZ}.
\begin{lemma}\label{lemma 2} {\bf (Lemma 3 in \cite{FNRZ})} Let $a\geq 0$.
\begin{enumerate}[i)]
 \item If $f,g\in{{\cal U}_a}$, then $\|fg\|_{{\cal U}_a}\leq C'\|f\|_{{\cal U}_a}\|g\|_{{\cal U}_a}$ .
 \item Let $T\in GL(d)$ with $|T|,|T|^{-1}\leq 2$. Then, for every $\delta>0$ and for every $f\in{\cal U}_a$, we have
 $\|\widehat{Tf}\|_{{\cal U}_{a-\delta}}\leq C'_{\delta}\|f\|_{{\cal U}_{a}}$.
 \item If $f\in {{\cal U}_a}$, then $\|\RR f\|_{{\cal U}_{a+d-2}}\leq C'\|f\|_{{\cal U}_a}$.
\end{enumerate}

\end{lemma}

\begin{lemma}\label{lemma 3}
{\bf (Lemma 5 in \cite{FNRZ})} Let $\beta>a$. For every $\sigma>0$, there exists a constant $C_0=C_{\sigma,a,\beta}>0$, such that
$\|f\|_{{\cal U}_a}\leq C_0\|f\|_{\infty}+\sigma\|f\|_{{\cal U}_{\beta}}$.
\end{lemma}
From Lemma \ref{lemma 2}, we get the following:
\begin{lemma}\label{lemma 4}
Let $a\geq 0$ and $f\in {{\cal U}_a}$. The following are true:
\begin{enumerate}[i)]
 \item $\|\Q f\|_{{\cal U}_{a+d-2}}\leq C''\|f\|_{{\cal U}_a}^{d-1}$.
 \item Let $T\in GL(d)$, with $|T|,|T|^{-1}\leq 2$. Then, for every $\delta>0$, we have $\|Tf\|_{{\cal U}_{a-\delta}}\leq C_{\delta}\|f\|_{{\cal U}_{a}}$.
\end{enumerate}

\end{lemma}
\beginproof (i) We notice that 
$$
\textnormal{det}_{(d-1)\times (d-1)}(x_1,\dots,x_{d-1})=\textnormal{det}_{d\times d}(x_1,\dots,x_{d-1}, \xi)
$$
for any $\xi \in \Sp^{d-1}$ and $x_1,\dots,x_{d-1} \in \xi^\perp$, thus
\begin{eqnarray*}\Q(g_1,\dots,g_{d-1})(\xi)&=&c_d\int\limits_{\Sp^{d-1}\cap \xi^{\bot}}\dots\int\limits_{\Sp^{d-1}\cap \xi^{\bot}}\det(x_1,\dots,x_{d-1})^2 \prod\limits_{k=1}^{d-1} g_k(x_k)dx_1\dots dx_{d-1}\\
&=&c_d\int\limits_{\Sp^{d-1}\cap \xi^{\bot}}\dots\int\limits_{\Sp^{d-1}\cap \xi^{\bot}}\textnormal{det}_{d\times d}(x_1,\dots,x_{d-1},\xi)^2 \prod\limits_{k=1}^{d-1} g_k(x_k)dx_1\dots dx_{d-1}\\
&=&c_d\int\limits_{\Sp^{d-1}\cap \xi^{\bot}}\dots\int\limits_{\Sp^{d-1}\cap \xi^{\bot}}\left(\sum\limits_{\sigma \in S_{d}}\sign(\sigma)  \xi_{\sigma(d)}\prod\limits_{k=1}^{d-1} x_{k,\sigma(k)}\right)^2 \prod\limits_{k=1}^{d-1} g_k(x_k)dx_1\dots dx_{d-1},
\end{eqnarray*}
where  $(x_{k,1}, \dots x_{k, d})$ are the coordinates of the vector $x_k\in \R^d$, $S_{d}$ is the set of permutations of $d$ elements and the last equality follows from the Leibniz formula for determinant.  Thus,

$\Q(g_1,\dots,g_{d-1})(\xi)$
\begin{eqnarray*}&=&c_d\int\limits_{\Sp^{d-1}\cap \xi^{\bot}}\!\!\!\!\dots\!\!\!\!\int\limits_{\Sp^{d-1}\cap \xi^{\bot}}\left(\sum\limits_{\sigma_1, \sigma_2 \in S_{d}}\sign(\sigma_1) \sign(\sigma_2)\xi_{\sigma_1(d)}\xi_{\sigma_2(d)} \prod\limits_{k=1}^{d-1} (x_{k,\sigma_1(k)}x_{k,\sigma_2(k)})\right) \prod\limits_{k=1}^{d-1} g_k(x_k)dx_1\dots dx_{d-1}\\
&=&\bar{c}_d\sum\limits_{\sigma_1, \sigma_2 \in S_{d}}\sign(\sigma_1) \sign(\sigma_2)\xi_{\sigma_1(d)}\xi_{\sigma_2(d)} \prod\limits_{k=1}^{d-1}  \int\limits_{\Sp^{d-1}\cap x^{\bot}}x_{k,\sigma_1(k)}x_{k,\sigma_2(k)}g_k(x_k)dx_k\\
&=&\bar{c}_d\sum\limits_{\sigma_1, \sigma_2 \in S_{d}}\sign(\sigma_1) \sign(\sigma_2)\xi_{\sigma_1(d)}\xi_{\sigma_2(d)} \prod\limits_{k=1}^{d-1}  \RR \left(x_{k,\sigma_1(k)}x_{k,\sigma_2(k)}g_k(x_k)\right)(\xi).
\end{eqnarray*}
Next we apply the triangular inequality to get that 
$$\|\Q f\|_{{\cal U}_{a+d-2}}\leq \bar{c}_d\sum\limits_{\sigma_1, \sigma_2 \in S_{d-1}}\left\|\xi_{\sigma_1(d)}\xi_{\sigma_2(d)} \prod\limits_{k=1}^{d-1}  \RR \left(x_{k,\sigma_1(k)}x_{k,\sigma_2(k)}g_k(x_k)\right)(\xi)\right\|_{{\cal U}_{a+d-2}}.
$$
 We use  Lemma \ref{lemma 2} (i) together with the fact that $\|\xi_i \xi_j\|_{a+d-2} \le c$ to get

$$\|\Q f\|_{{\cal U}_{a+d-2}}\leq \overline{C_0}\bigg(\max_{1\le i,j \le d}\|\RR(x_{i}x_{j}f(x))\|_{{\cal U}_{a+d-2}}\bigg)^{d-1} .$$
Now, by Lemma \ref{lemma 2} (iii) we derive:
$$\|\Q f\|_{{\cal U}_{a+d-2}}\leq \overline{C_0}C'^{d-1}\bigg(\max_{1\le i,j \le d}\|x_{i}x_{j}f(x)\|_{{\cal U}_{a}}\bigg)^{d-1}  $$
and again by Lemma \ref{lemma 2} (i), we conclude:
$$\|\Q f\|_{{\cal U}_{a+d-2}}\leq \overline{C_0}C'^{d-1}C'^{d-1}\bigg(\max_{1\le i,j \le d}\|x_{i}x_{j}\|_{{\cal U}_{a}}\bigg)^{d-1}
\|f\|_{{\cal U}_{a}}^{d-1}=
C''\|f\|_{{\cal U}_{a}}^{d-1}.$$
(ii) Fix $\delta>0$. By Lemma \ref{lemma 2} (ii), we have $\|\widehat{Tf}\|_{{\cal U}_{a-\delta}}\leq \|f\|_{{\cal U}_{a}}$, so
Lemma \ref{lemma 2} (i) implies:
$$\|Tf\|_{{\cal U}_{a-\delta}}\leq C'\||Tx|^{-d}\|_{{\cal U}_{a-\delta}}\|f\|_{{\cal U}_{a-\delta}}\leq
C'C'_{\delta}\||Tx|^{-d}\|_{{\cal U}_{a-\delta}}\|f\|_{{\cal U}_{a}}.$$
For $f\equiv 1$, we have $\widehat{Tf}=|Tx|^{-1}$, thus again by Lemma \ref{lemma 2} (i), (ii), we get:
$$\||Tx|^{-d}\|_{{\cal U}_{a-\delta}}\leq (C')^d\||Tx|^{-1}\|_{{\cal U}_{a-\delta}}^d\leq (C')^d(C'_{\delta})^d\|1\|_{{\cal U}_{a}}^d=(C')^d(C'_{\delta})^d, $$
which finishes the proof.
\endproof
\begin{lemma}\label{lemma 6}
Fix $d>2$ and let $a, L, C$ be the constants from Lemma \ref{lemma 1}. There exist $0<\lambda<1$ and $\epsilon>0$, with the following property:
If $\varepsilon\in (0,\epsilon)$ and $f:\Sp^{d-1}\rightarrow \mathbb{R}$ is a function with $f=1+\varphi$, $\int_{\Sp^{d-1}}\varphi=0$,
$\|\varphi\|_2<\varepsilon$ and $\|\varphi\|_{{\cal U}_{a}}<1/L$, then there exists $T\in GL(d)$ and $\gamma>0$, such that the function
$\widetilde{f}=\gamma \Q(Tf)$ can be written as $\widetilde{f}=1+\widetilde{\varphi}$, for some $\widetilde{\varphi}$ with $\int_{\Sp^{d-1}}\widetilde{\varphi}=0$, $\|\widetilde{\varphi}\|_{{\cal U}_a}<1/L$ and $\|\widetilde{\varphi}\|_2<\lambda \varepsilon$.
\end{lemma}
\beginproof Use Lemma \ref{lemma 1} to find $T\in GL(d)$, such that $Tf=1+\psi$, $\|\psi\|_2\leq \varepsilon+C\varepsilon^{\frac{d+5}{d+3}}$, $\|\psi\|_{\infty}\leq C\varepsilon^{\frac{4}{d+3}}$ and
$\|H_2^{\psi}\|_2\leq C\varepsilon^{\frac{d+5}{d+3}}$.

Using rotation invariance of Haar measure we get that 
$$\int_{\Sp^{d-1}\cap \xi^\perp}\dots\int_{\Sp^{d-1}\cap \xi^\bot }|\textnormal{det}(x_1\dots,x_{d-2},x)|^2 dx_1\dots dx_{d-2}$$ 
is a constant for all $x \in \Sp^{d-1}\cap \xi^\perp$.  Thus, it follows that for any function $g$, we have:
\begin{equation}
\Q(g,1[d-2])=\RR(g).\label{observation}
\end{equation}
One should also notice, that the constant in the above equality is exactly 1, by setting $g=1$. Therefore
\begin{eqnarray*}
\Q(Tf)&=&\Q(1+\psi)\\
&=&\sum_{k=0}^{d-1}\binom{d-1}{k}\Q(\psi[k],1[d-1-k])\\
&=&1+(d-1)\RR(\psi)+\sum_{k=2}^{d-1}\binom{d-1}{k}\Q(\psi[k],1[d-1-k]).
\end{eqnarray*}
Let $\nu=\sum_{k=2}^{d-1}\binom{d-1}{k}\Q(\psi[k],1[d-1-k])$. Then, 
\begin{eqnarray}
\|\nu-H_0^{\nu}\|_2\leq\|\nu\|_2&\leq&\sum_{k=2}^{d-1}\binom{d-1}{k}\|\Q(|\psi|\ [k],1[d-1-k])\|_2\nonumber\\
&\leq& \sum_{k=2}^{d-1}\binom{d-1}{k}\|\psi\|_{\infty}^{k-1}\|\Q(|\psi|,1[d-2])\|_2  \nonumber\\
&=&\sum_{k=2}^{d-1}\binom{d-1}{k}\|\psi\|_{\infty}^{k-1}\|\RR(|\psi|)\|_2\nonumber\\
&\leq& \sum_{k=2}^{d-1}\binom{d-1}{k}\|\psi\|_{\infty}^{k-1}\|\psi\|_2\nonumber\\
&\leq& (d-1)\max_k\binom{d-1}{k}(1+C)^d\varepsilon^{\frac{4}{d+3}}C\varepsilon\nonumber\\
&=&C_1\varepsilon^{\frac{d+7}{d+3}}\label{proof-eq-1}.
\end{eqnarray}
Similarly, we have:
\begin{eqnarray}
\|\nu-H_0^{\nu}\|_{\infty}\leq 2\|\nu\|_{\infty}&\leq& 2\sum_{k=2}^{d-1}\binom{d-1}{k}\|Q(|\psi|\ [k],1[d-1-k])\|_{\infty}\nonumber\\
&\leq& 2\sum_{k=2}^{d-1}\binom{d-1}{k}\|\psi\|_{\infty}^k\leq C_2\varepsilon^{\frac{8}{d+3}}.
\label{proof-eq-2}
\end{eqnarray}
Write 
$$
\Q(Tf)=1+(d-1)\RR(\psi)+\nu=1+(d-1)H_0^{\psi}+H_0^{\nu}+(d-1)\RR H_2^{\psi}+(d-1)\RR(\psi-H_0^{\psi}-H_2^{\psi})+(\nu-H_0^{\nu})
$$
and set $\zeta=(d-1)H_0^{\psi}+H_0^{\nu}$. Since $\zeta$ is a constant, we have:
$$|\zeta|\leq (d-1)|H_0^{\psi}|+|H_0^{\nu}|\leq (d-1)\|\psi\|_2+\|\nu\|_2\leq C_3\varepsilon.$$
Also, $(d-1)\|\RR H_2^{\psi}\|_2\leq (d-1)\|H_2^{\psi}\|_2\leq C \varepsilon^{\frac{d+5}{d+3}}$ and since decomposition of $\RR(\psi)$ into spherical harmonics does not contain  spherical harmonics of odd degree, we have: 
\begin{eqnarray}
(d-1)\|\RR(\psi-H_0^{\psi}-H_2^{\psi})\|_2
&=&(d-1)\sum_{k=4}^{\infty}v_{k,d}\|H_k^{\psi}\|_2\nonumber \\ &\leq& (d-1)v_{4,d}\sum_{k=4}^{\infty}\|H_k^{\psi}\|_2\nonumber\\
&=&\frac{3}{d+1}\|\psi-H_0^{\psi}-H_2^{\psi}\|_2\nonumber \\&\leq& \frac{3}{4}\|\psi\|_2\leq\frac{3}{4}(\varepsilon+C\varepsilon^{\frac{d+5}{d+3}})\nonumber\\ &\leq& C'\varepsilon,\label{proof-eq-3}
\end{eqnarray}
provided that $\varepsilon$ is small enough. Take $\gamma=(1+\zeta)^{-1}=O(\varepsilon)$ and set
$$
\widetilde{\varphi}:=\gamma\big[(d-1)RH_2^{\psi}+(d-1)\RR(\psi-H_0^{\psi}-H_2^{\psi})+(\nu-H_0^{\nu})\big].$$
Then, by (\ref{proof-eq-1}), (\ref{proof-eq-3}) and the fact that $(d-1)\|\RR H_2^{\psi}\|_2\leq C \varepsilon^{\frac{d+5}{d+3}}$, we immediately get:
$$\|\widetilde{\varphi}\|_2\leq (1+C_2\varepsilon)\big(C \varepsilon^{\frac{d+5}{d+3}}+C'\varepsilon+C_1\varepsilon^{\frac{d+7}{d+3}}\big)\leq C\varepsilon,$$
provided that $\varepsilon$ is sufficiently small. Moreover, by (\ref{proof-eq-2}), we get:
$$\|\widetilde{\varphi}\|_{\infty}\leq \gamma\big((d-1)\|\RR(\psi-H_0^{\psi})\|_{\infty}+\|\nu-H_0^{\nu}\|_{\infty}\big)
\leq \gamma\big((d-1)2\|\psi\|_{\infty}+C_2\varepsilon^{\frac{4}{d+3}})\leq C_4\varepsilon^{\frac{4}{d+3}}.$$
Finally, to estimate $\|\widetilde{\varphi}\|_{{\cal U}_a}$, note that $|T|,|T|^{-1},\|f\|_{{\cal U}_a}<2$. Fix any $0<\delta<1$.
Then, by Lemma \ref{lemma 4} we get $\|Tf\|_{{\cal U}_{a-\delta+d-2}}\leq \|Tf\|_{{\cal U}_{a-\delta}}\leq C_{\delta}\|f\|_{{\cal U}_a}<2C_{\delta}$. Again, Lemma \ref{lemma 4} implies 
$$
\|\gamma Q(Tf)\|_{{\cal U}_{a-\delta+d-2}}\leq \gamma 2 C_{\delta}C''\|f\|_{{\cal U}_a}^{d-1}\leq \gamma 2 C_{\delta}C''2^{d-1}=:C_5.
$$
 Therefore, 
 $$
 \|\widetilde{\varphi}\|_{{\cal U}_{a-\delta+d-2}}=\|\gamma Q(Tf)-1\|_{{\cal U}_{a-\delta+d-2}}\leq C_5+1=:C_6.
 $$
  Now, one may use Lemma \ref{lemma 3} with $\sigma=1/(2C_6L)$, to conclude:
$$\|\widetilde{\varphi}\|_{{\cal U}_a}\leq C_{a,a-\delta+d-2}\|\widetilde{\varphi}\|_{\infty}+\sigma\|\widetilde{\varphi}\|_{{\cal U}_{a+\delta+d-2}}\leq C_7\varepsilon^{\frac{4}{d+3}}+C_6/(2C_6L)\leq 1/L,$$
provided that $\varepsilon$ is small enough. 
\endproof
\noindent{\bf  Proof of Theorem \ref{main-theorem}:}  
Since the ``iteration'' Lemma \ref{lemma 6} is established, the rest of the proof of Theorem \ref{main-theorem} is exactly the same as the proof of Theorem A. For the sake of completeness, we will briefly repeat the argument. Let $K$ be a convex body with absolutely continuous surface area measure and $\|f_K-1\|_{\infty}<\varepsilon$, for some $0<\varepsilon<1$. Then, $\|f_K\|_{{\cal U}_0}<1+\varepsilon$ and by Lemma \ref{lemma 4} we easily get
$\|f_{\Pi^mK}\|_{{\cal U}_{m(d-2)}}<C^m(1+\varepsilon)^m$. Choosing $m$  such that $a<m(d-2)<2a$, where $a$ is the constant from Lemma \ref{lemma 6}, we conclude that $\|f_{\Pi^mK}\|_{{\cal U}_{\beta}}<\overline{C}$, where $\beta=m(d-2)$. On the other hand, it is clear from the definition of $\Q(f)$ that since $1-\varepsilon\leq f_{K}\leq 1+\varepsilon$, we have: 
$$
(1-\varepsilon)^{(d-1)^m}\leq f_{\Pi^mK}\leq (1+\varepsilon)^{(d-1)^m}.
$$
 So, $f_{\Pi^mK}/(\int_{\Sp^{d-1}}f_{\Pi^mK})=1+\varphi$, with $\int_{\Sp^{d-1}}\varphi =0$ and 
 $$
 \|\varphi\|_{{\cal U}_{\beta}}\leq 1+\|f_{\Pi^mK}\|_{{\cal U}_{\beta}}/(\int_{\Sp^{d-1}}f_{\Pi^mK})<\overline{C}'.
 $$ 
 Clearly, $\|\varphi\|_{\infty}\leq \overline{C}''\varepsilon$, so by Lemma \ref{lemma 3}, we get: $\|\varphi\|_{{\cal U}_a}\leq C_{\sigma,a,\beta}\overline{C}''\varepsilon+\sigma\overline{C}'\leq 1/L$, provided that $\varepsilon$ is small enough. Therefore, replacing $K$ with $\Pi^mK$, we may assume that the assumptions of Lemma \ref{lemma 6} hold. Since the factor $|\textnormal{det} T|^{(d-1)/d}$ can be eliminated by rescaling, Lemma \ref{lemma 6} implies that for each $m\in \mathbb{N}$, there exists $T_m\in GL(d)$, such that $\|f_{T_m(\Pi^mK)}-1\|_2\rightarrow 0$, as $m\rightarrow \infty$ and $\|f_{T_m(\Pi^mK)}\|_{{\cal U}_a}<1/L$. Now, Lemma \ref{lemma 1} (i) implies that $\|f_{T_m(\Pi^mK)}-1\|_{\infty}\rightarrow 0$, as required.$\Box$
\section{Proof of Theorem \ref{corollary 2}}

In this section, some facts from the theory of mixed volumes and quermassintegrals of convex bodies will be used. We refer to \cite{Sch1} for an extensive discussion on this topic. Let $m\in\mathbb{N}$, $t_1,\dots,t_m\in\mathbb{R}$ and $K_1,\dots,K_m$ be convex bodies in $\mathbb{R}^d$. Recall the classical Minkowski formula about mixed volumes:
\begin{equation}
|t_1K_1+\dots t_mK_m|=\sum_{i_1,\dots,i_d=1}^mt_{i_1}\dots t_{i_d}V(K_{i_1},\dots,K_{i_d}).
\label{mixed-vol} 
\end{equation}
Note, also, that $W_j(K)=V(K[d-j],B_2^d[j])$ ($K$ appears $d-j$ times and $B_2^d$ appears $j$ times), $j=0,\dots,d-1$.
For instance, $W_0(K)$, $W_1(K)$ and $W_{d-1}(K)$ are proportional to the volume, the surface area and the mean width of $K$ respectively.

We will establish an inequality of the form $$W_{d-2}(\Pi_{d-2}K)\geq a_dW^2_{d-1}(K)+\beta_dW_{d-2}(K),$$ where $a_d$ and $\beta_d$ are positive constants that depend only on $d$ and we will prove that there is equality when $K=B_2^d$. The computation of the exact values of $a_d$ and $\beta_d$ is an elementary (but rather tedious) task and it is left to the reader. In what follows, every constant that appears will depend only on the dimension. 


 
In order to prove Theorem \ref{corollary 2}, we would like to be able to claim that  functions 
$$
t\mapsto |B_2^d+tK| \mbox{ and  }t\mapsto|\Pi (B_2^d+tK)|
$$ 
well defined and  twice differentiable at $t=0$. Unfortunately, this cannot be expected to be true for a general convex body $K$. To overcome this difficulty, we will first prove our theorem for symmetric convex bodies $K$ with smooth boundary. It is  well known (see \cite{W3} or again \cite{Sch1}) that in this case  $K$ is a generalized zonoid, i.e. the  support function of $K$ is given by
$$h_K(u)=\frac{1}{2}\int_{\Sp^{d-1}}|\langle x,u \rangle|f(x)dx,$$
where $f:\Sp^{d-1} \to \R$ is an even, continuous function. 

Let $\delta = (\max\limits_{x\in\Sp^{d-1}} |f(x)| )^{-1}$. For  each $t\in (-\delta, \delta)$ we   define a convex body $M_t$ via its support function
\begin{equation}
h_{M_t}(u)=\frac{1}{2}\int_{\Sp^{d-1}}|\langle x,u\rangle|(1+tf(x))dx, \qquad u\in \Sp^{d-1}.\label{equation-M_t}
\end{equation}
Note that for $t>0$, $M_t=\Pi B_2^d+tK$, while for $t<0$ Minkowski addition is replaced with Minkowski subtraction. Although Minkowski subtraction, when meaningful, does not always have good behavior (see \cite{Sch1}), the integral representation (\ref{equation-M_t}) (since $1+tf\geq 0$ for $|t|<\delta$), allows us, by a result of Weil, to derive nice integral expressions for the volumes of $M_t$ and $\Pi M_t$ (see \cite{W1} or \cite{Sch1} for $|M_t|$ and \cite{Sa1} for $|\Pi M_t|$). Those formulas give us  that the functions $|M_t|$ and $|\Pi M_t|$ are polynomials in $t$, hence both twice differentiable in $(-\delta,\delta)$. This simple observation will be crucial for the proof of Theorem \ref{corollary 2}. 

Moreover, we also note that  the curvature function  $f_{M_t}$ is proportional to $\Q(1+tf)$. Thus,   there exists $\delta'>0$ such that     $f_{M_t}$ is sufficiently close to $1$ for  $|t|<\delta'$. Next, applying   Theorem \ref{corollary 1}, we conclude that $P(M_t)\geq P(B_2^d)$ for $|t|<\delta'$.
%
Thereby, it follows by an obvious rescaling that the point $t=0$ is a local minimizer for the function $t\mapsto P(B_2^d+tK)$. One can check that if for two twice differentiable functions $f,g:\mathbb{R}\to \mathbb{R_+}$ with $g(0),g''(0)>0$, the function $f/g$
is attains a local minimum at $0$, then $f''(0)/g''(0)\geq f(0)/g(0)$ and there is equality if $f/g$ is constant. Therefore,
\begin{equation}
\frac{d^2}{dt^2}|\Pi(B_2^d+tK)|\bigg|_{t=0}\geq \frac{|\Pi B_2^d|}{|B_2^d|^{d-1}}\frac{d^2}{dt^2}\big(|B_2^d+tK|^{d-1}\big)\bigg|_{t=0}=:A_d\frac{d^2}{dt^2}\big(|B_2^d+tK|^{d-1}\big)\bigg|_{t=0},\label{eq1-proof-of-corollary-2}
\end{equation}
with equality if $K=B_2^d$. By definition, the first and second derivative of the function $t\mapsto |B_2^d+tK|$ at $t=0$ are proportional to $W_{d-1}(K)$ and $W_{d-2}(K)$ respectively. Hence an elementary calculation shows that
$$\frac{d^2}{dt^2}\big(|B_2^d+tK|^{d-1}\big)\bigg|_{t=0}=\gamma_dW_{d-1}^2(K)+\delta_dW_{d-2}(K).$$
Next, note that for $u\in \Sp^{d-1}$, using (\ref{mixed-vol}), we have:
\begin{eqnarray*}
h_{\Pi (B_2^d+tK)}(u)=|(B_2^d+tK)|u^{\bot}|&=&|t(K|u^{\bot})+B_2^{d-1}|\\
&=&|B_2^{d}|u^{\bot}|+\sum_{i=1}^{d-2}\zeta_{i,d}t^iW_{d-i-1}(K|u^{\bot})\\
&=&\sum_{i=1}^{d-2}\zeta_{i,d}t^ih_{\Pi_{d-i-1}K}(u)+h_{\Pi B_2^d}(u),\qquad t\geq 0,
\end{eqnarray*}
where $\zeta_{i,d}$ are  positive constants. Therefore, again by (\ref{mixed-vol}), we get:
\begin{eqnarray*}
|\Pi(B_2^d+tK)|&=&\Big|\sum_{i=1}^{d-2}\zeta_{i,d}t^i\Pi_{d-i-1}K+\zeta_{0,d} B_2^d\Big|\\
&=&t^2\Big(\eta_dW_{d-2}(\Pi_{d-2}K)+\theta_dW_{d-1}(\Pi_{d-3}K)\Big)+\Theta(t)+O(t^3)\\
&=&t^2\Big(\eta_dW_{d-2}(\Pi_{d-2}K)+\lambda_dW_{d-2}(K)\Big)+\Theta(t)+O(t^3),\qquad t\geq 0,
\end{eqnarray*}
where we used (\ref{identity}). Consequently,
$$\frac{d^2}{dt^2}\big|\Pi(B_2^d+tK)\big|\bigg|_{t=0}=2\eta_dW_{d-2}(\Pi_{d-2}K)+2\lambda_dW_{d-2}(K),$$
so (\ref{eq1-proof-of-corollary-2}) becomes:
$$W_{d-2}(\Pi_{d-2}K)\geq \frac{\gamma_d}{2\eta_d}W^2_{d-1}(K)+\frac{\delta_d-2\lambda_d}{2\eta_d}W_{d-2}(K)=:a_dW^2_{d-1}(K)+\beta_dW_{d-2}(K).$$
Note that since (\ref{eq1-proof-of-corollary-2}) is sharp for the ball, the last inequality is also sharp. To conclude that $\beta_d>0$, notice that if $\beta_d\leq 0$, then by the use of (\ref{A-F-inequalities}) and (\ref{identity}), for $i=d-2$, we would get:
$$a_dW_{d-1}^2(K)\leq  W_{d-2}(\Pi_{d-2}K)+|\beta_d|W_{d-2}(K)\leq \omega_dW_{d-1}^2(\Pi_{d-2}K)+|\beta_d|\omega_dW^2_{d-1}(K)
=:a'_dW_{d-1}^2(K),
$$
for all convex bodies $K$, with equality everywhere if $K$ is a ball. This is a contradiction, thus $\beta_d>0$.

We have established our inequality for symmetric convex bodies with smooth boundary. This extends easily to all symmetric convex bodies. Assume, now, that $K$ is a general convex body. Define the difference body $DK=\frac{1}{2}(K-K)$. Since $DK$ is symmetric, we have:$$W_{d-2}(\Pi_{d-2}DK)\geq a_dW^2_{d-1}(DK)+\beta_dW_{d-2}(DK).$$
By the additivity of 
the mean width functional, 
we immediately conclude that $\Pi_{d-2}DK=\Pi_{d-2}K$ and $W_{d-1}(DK)=W_{d-1}(K)$. Hence, we only need to prove that $W_{d-2}(DK)\geq W_{d-2}(K)$ and the inequality is strict if $K$ is not centrally symmetric.
By 
the additivity of mixed volumes, we easily get:
\begin{eqnarray*}
W_{d-2}(DK)=\frac{1}{2}V(K,K,B_2^d[d-2])+\frac{1}{2} V(K,-K,B_2^d[d-2]).
\end{eqnarray*}
Finally, use the general Aleksandrov-Fenchel inequality $$V(K_1,\dots,K_d)^2\geq V(K_1,K_1,K_3,\dots,K_d)V(K_2,K_2,K_3,\dots,K_d),$$
to conclude that
$$V(K,-K,B_2^d,\dots, B_2^d)\geq V(K,K,B_2^d,\dots, B_2^d)^{\frac{1}{2}}V(-K,-K,B_2^d,\dots, B_2^d)^{\frac{1}{2}}=W_{d-2}(K).$$
It is well known (see \cite[Theorem 7.6.2]{Sch1}) that there is equality in the last inequality if and only if $K$ is centrally symmetric, which ends our proof. $\Box$
\\
\\
Before ending this note, we would like to remark that our method does not provide further information about the equality cases for the inequality of Theorem \ref{corollary 2}. One may naturally conjecture that there is equality if and only if $K$ is a ball.

\vspace{0.8cm}

\noindent Christos Saroglou \\
Department of Mathematical Sciences \\
Kent State University \\
Kent, OH 44242, USA \\
E-mail address: csaroglo@kent.edu \ \&\ christos.saroglou@gmail.com

\vspace{0.8cm}

\noindent Artem Zvavitch \\
Department of Mathematical Sciences \\
Kent State University \\
Kent, OH 44242, USA \\
E-mail address: zvavitch@math.kent.edu


\begin{thebibliography}{999}

\bibitem[FNRZ]{FNRZ}{\sc A. Fish, F. Nazarov, D. Ryabogin, A. Zvavitch,} {\em The behavior of iterations of the intersection body operator in a small neighborhood of the unit ball}, Adv. Math., 226 (2011), no 3, 2629-2642.

\bibitem[Ga]{G}{\sc R. J. Gardner}, {\em Geometric tomography}. Second edition,
Encyclopedia of Mathematics and its Applications 58 Cambridge University Press, Cambridge, 2006.

\bibitem[Gr]{Gr} {\sc H. Groemer,} {\em 
Geometric Applications of Fourier Series and Spherical Harmonics,} Cambridge University Press, New York, 1996.

\bibitem[GP]{GP}{\sc A. Giannopoulos, M. Papadimitrakis}, {\em Isotropic surface area measures}, Mathematika 46 (1999) 1-�13.


\bibitem[Ko]{Ko} {\sc A.~Koldobsky},
{\em Fourier Analysis in Convex Geometry}, Math. Surveys and Monographs,
AMS, Providence RI 2005.



\bibitem[KoY]{KoY}{\sc A.~Koldobsky, V.~Yaskin},  {\em The Interface between Convex Geometry and Harmonic Analysis}, CBMS Regional Conference Series, 108, American Mathematical Society, Providence RI, 2008.


\bibitem[Lu1]{Lu1} {\sc E.~Lutwak}, {\em Intersection bodies and dual mixed volumes},
Advances in Math. 71 (1988), 232--261.


\bibitem[Lu2]{Lut3}{\sc E. Lutwak}, {\em Centroid bodies and dual mixed volumes}, Proc. London Math. Soc. 60 (1990), 365--391.

\bibitem[Lu3]{Lut2}{\sc E. Lutwak}, {\em On quermassintegrals of mixed projection bodies},
Geom. Dedicata 33 (1990) 51--58.

\bibitem[Lu4]{Lu4} {\sc E.~Lutwak}, {\em Selected affine isoperimetric inequalities}, In a Handbook of Convex Geometry, ed. by P.M. Gruber and J.M. Wills. North-Holland, Amsterdam, 1993, pp. 151-176.


\bibitem[P1]{PE2}{\sc C. M. Petty}, {\em Projection bodies}. 1967 Proc. Colloquium on Convexity (Copenhagen, 1965) 234--241 Kobenhavns Univ. Mat. Inst., Copenhagen.
\bibitem[P2]{PE}{\sc C. M. Petty}, {\em Isoperimetric problems.} Proceedings of the Conference on Convexity and Combinatorial Geometry (Univ. Oklahoma, Norman, Okla., 1971) 26--41. Dept. Math., Univ. Oklahoma, Norman, Okla., 1971.
\bibitem[Sa]{Sa1}{\sc C. Saroglou}, {\em Volumes of projection bodies of some classes of convex bodies,} Mathematika 57 (2011) 329--353.

\bibitem[Sch1]{Sch3}{\sc R. Schneider}, {\em Geometric inequalities for Poisson processes of convex bodies and
cylinders}, Results Math. 11 (1987) 165-185.
\bibitem[Sc2]{Sch1} {\sc R. Schneider,} {\em Convex bodies: the Brunn-Minkowski theory}. Second expanded edition.
Encyclopedia of Mathematics and its Applications 151, Cambridge University Press, Cambridge, 2014.
\bibitem[W1]{W1}{\sc W. Weil}, {\em $\ddot{\textnormal{U}}$ber die
Projektionenk$\ddot{\textnormal{o}}$rper konvexer Polytope},
Arch. Math. 22 (1971) 664-672.
\bibitem[W2]{W2}{\sc W. Weil}, {\em Kontinuierliche Linearkombination won
Strecken}, Math. Z. 148 (1976) 71-84.
\bibitem[W3]{W3}{\sc W. Weil}, {\em Centrally  symmetric  convex  bodies  and  distributions}, Israel  J.  Math. 24 (1976), 352–--367.
\bibitem[Zh]{Zh}{\sc G. Zhang}, {\em The affine Sobolev inequality}, J. diff. geom. 53 (1999), 183-202.



\end{thebibliography}
\end{document}